\documentclass[12pt]{article}
\usepackage{amsmath,amssymb,amsthm}
\usepackage{graphicx}
\input epsf
\input epsf

 \def\R{{\mathbb R}}  
\long\def\comment#1\endcomment{}

\begin{document}


\centerline{\uppercase{\bf  On intersection of three embedded
spheres in 3-space}}
 \bigskip
\centerline{\bf A. Rukhovich}


\bigskip
\small
{\bf Abstract.} We study intersection of two polyhedral spheres without self-intersections in 3-space $\R^3$. We find necessary and sufficient conditions on sequences $\vec x = (x_1,x_2,\dots,x_n)$, $\vec y = (y_1,y_2,\dots, y_n)$ of positive integers, for existence of polyhedral spheres  $f,g\subset\R^3$ such that

$\bullet$ $f-g$ has $n$ connected components, which can be numbered so that the $i$-th component has $x_i$ neighbors in $f$ and

$\bullet$ $g-f$ has $n$ connected components, which can be numbered so that the $i$-th component has $y_i$ neighbors in $g$.

Analogously we study intersection of {\it three} polyhedral spheres without self-intersections in 3-space.

\comment

We study intersection of two polyhedral spheres without self-intersections in 3-space. We find necessary and sufficient conditions on sequences  x = x_1,x_2,...,x_n,  y = y_1,y_2,...,y_n  of positive integers, for existence of   polyhedral spheres  f,g in R^3   such that

* f-g has n connected components, which can be numbered so that the $i$-th component has x_i neighbors in f and

* g-f has n connected components, which can be numbered so that the $i$-th component has y_i neighbors in g.

Analogously we study intersection of  three  polyhedral spheres without self-intersections in 3-space.

\endcomment

\normalsize
\bigskip
{\bf 1. Introduction and main results.}

\smallskip
{\bf Theorem 1.} {\it Let $n$ be a positive integer and $\vec x = (x_1,x_2,\dots,x_n)$, $\vec y = (y_1,y_2, \dots, y_n)$ be sequences of
positive integers.
There exist
  polyhedral spheres $f,g\subset\R^3$  intersecting by closed broken lines and such that

$\bullet$
$f-g$ has $n$ connected components, which can be numbered so that the $i$-th component has $x_i$ neighbors in $f$;

$\bullet$
$g-f$ has $n$ connected components, which can be numbered so that the $i$-th component has $y_i$ neighbors in $g$;

if and only if $\sum_{i=1}^nx_i=\sum_{i=1}^ny_i=2n-2$.}

\smallskip

Two connected components of $f-g$ are {\it neighbors in $f$} if their closures intersect.
Analogously one defines {\it neighbors in $g$} for components of $g-f$.

In this note a {\it sphere} is a polyhedral sphere, i.e. a polyhedron homeomorphic to the
sphere. A {\it circle} is a closed simple broken line in a sphere.

\begin{figure}[h]
  \centerline{\includegraphics[width=8cm]{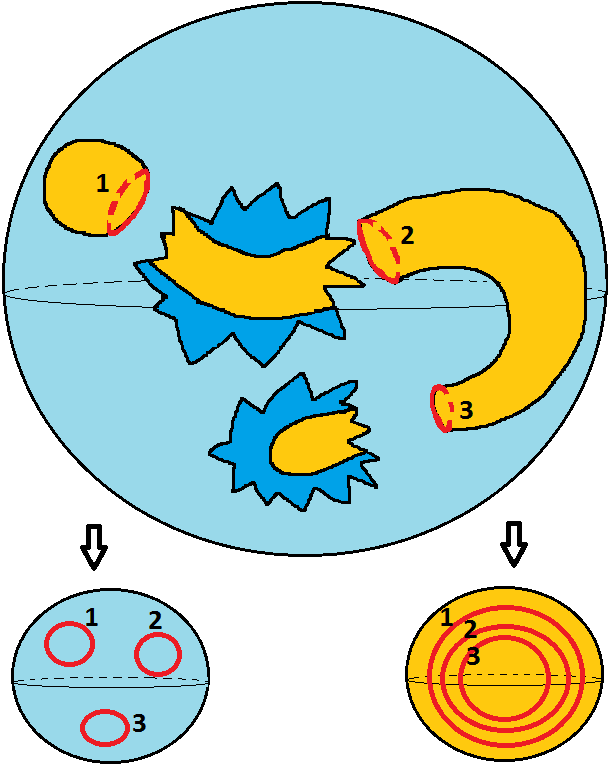}}
\caption{Intersecting spheres}\label{example}
\end{figure}

\smallskip
{\bf Remark.}
(1) E.g. in Figure \ref{example} the intersection of  spheres $f, g$ consists of three circles that split each sphere into four components.
One component of $f-g$ has three neighbors, each of the other components  has one neighbor.
Two components of $g-f$ have each two neighbors, each of the other two components has one neighbor.
For an elementary exposition of this and related results see [T], in russian [Tr].

\begin{figure}[h]
  \centerline{\includegraphics[width=10cm]{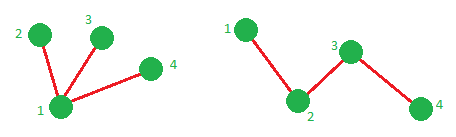}}
\caption{Graphs corresponding to Figure \ref{example}}\label{graphs}
\end{figure}

(2) Let $f,g\subset\R^3$ be  spheres.
Define graph $F$ as follows. The vertices are connected components of $f-g$. Two vertices are connected
by an edge if the corresponding connected components are neighbors. The graph $F$ may be called {\it dual} to the family of circles $f\cap g$ in $f$. Analogously define graph $G$. Then Theorem 1 describes pairs of degree sequences of such graphs.
In Figure \ref{graphs} there are graphs $F, G$ for  spheres from Figure \ref{example}.

(3) The `only if' part of Theorem 1 is essentially known and is essentially proved in [N]
(we present elementary proofs).
The `if' part is new.
The `if' part is interesting because there exist two unions of circles on spheres which are not realizable by
intersecting spheres in $\R^3$ [A].

(4) The conditions in Theorem 1 can be reformulated as follows:


$\bullet$ $f-g$ is the disjoint union of a sphere
with $x_1$ holes, a sphere with $x_2$ holes, \dots, a sphere with $x_n$ holes;

$\bullet$ $g-f$ is the disjoint union of a sphere with $y_1$
holes, a sphere with $y_2$ holes, \dots, a sphere with $y_n$ holes.

\smallskip
We call a sequence of positive integers $x_1, x_2, \dots, x_n$ {\it tree-like} if $\sum_{i=1}^nx_i=2n-2$.

\smallskip
{\bf Theorem 2.} {\it Let $n_1, n_2, n_3$ be positive integers and $$x_{1 1}, x_{1 2}, \dots, x_{1 n_1},\quad x_{2 1}, x_{2 2}, \dots ,x_{2 n_2},\quad x_{3 1}, x_{3 2}, \dots, x_{3 n_3}$$ be sequences of positive integers. There exist   spheres $f_1, f_2, f_3\subset\R^3$  pairwise intersecting by closed broken lines and  such that

$\bullet$ $f_1 \cap f_2 \cap f_3 =\varnothing$;

$\bullet$ $f_k-f_{k+1}-f_{k+2}$ has $n_k$ connected components, which can be numbered so that the $i$-th component has $x_{k i}$ neighbors in $f_k$, for each $1\le k \le 3$

if and only if the sequences are tree-like, $n_1+n_2+n_3$ is odd and $n_k < n_{k+1}+n_{k+2}$ for each $1\le k \le 3$.}

\smallskip
Here subscripts $k, k+1, k+2$ are considered mod 3.

\smallskip
{\bf Remark.}
(2) Theorem 2 has a graph-theoretic interpretation just as Theorem 1. 

(3) The `only if' part of Theorem 2 is trivial. The `if' part is not trivial and new. The `if' part is proved using Theorem 1.

(4) The second condition in Theorem 2 can be reformulated as follows:

$\bullet$ $f_k-f_{k+1}-f_{k+2}$ is the disjoint union of a sphere with $x_{k 1}$ holes, a sphere with  $x_{k 2}$ holes, \dots, a sphere with $x_{k n_k}$ holes for each $1\le k \le 3$.


\smallskip
{\it Proof of the `only if' part in Theorem 1.}
Recall definition of a graph $F$.
The vertices of $F$ are connected components of $f-g$.
Two vertices are connected
by an edge if the corresponding connected components are neighbors. Denote by $n$ the number of the vertices. The number of
the edges is equal to the number of circles in $f\cap g$. This number is $\sum_{i=1}^nx_i/2$. It is obvious that
$F$ is connected. By the Jordan Curve Theorem, $F$ is split by any
vertex.
So $F$ is a tree. Hence the number of edges is
$n-1=\sum_{i=1}^nx_i/2$. QED

 \smallskip
{\it Proof of the `only if' part in Theorem 1 suggested by T. Nowik.}
By induction on the number of circles. The statement
is true for one circle (there are only 2 disks on each sphere hence $n=2$). Each
additional circle splits one component into two, and adds two boundary circles.
QED

 \smallskip
{\it Proof of the `only if' part of Theorem 2.}
The necessity of the first condition obviously follows from the `only if' part of Theorem 1.

Let $m_3, m_2, m_1$ be the numbers of the circles in $f_1\cap f_2$, $f_1\cap f_3$ and $f_2\cap f_3$. Then $n_1=m_3+m_2+1$, $n_2=m_3+m_1+1$, $n_3=m_2+m_1+1$.

So $n_1+n_2+n_3=2(m_3+m_2+m_1)+3$ is odd.

Since $2m_k+1>0$ we have $n_k < n_{k+1}+n_{k+2}$ for each $1\le k\le 3$.
QED

\bigskip
{\bf 2. Proofs of the `if' parts of Theorems 1 and 2.}

\smallskip

A pair $(\vec x, \vec y)$ of sequences of positive integers is called {\it strongly realizable} if there exist two curved spheres $S,T$

(1) whose intersection consists of $n-1$ circles and splits

\quad $\bullet$ $S$ into $n$ connected components which can be numbered so that the $i$-th connected component has $x_i$ neighbors in $S$, and

\quad $\bullet$ $T$ into $n$ connected components which can be numbered so that the $i$-th connected component has $y_i$ neighbors in $T$;

(2) there is a circle of $S\cap T$ that bounds a disk and a component with $x_1$ neighbors in $S-T$,
as well as bounds a disk and a  component with $y_1$ neighbors in $T-S$.

Note that if $x_1=1$, then `the disk and the component with $x_1$ neighbors' in (2) could be the same component.
Same remark can be done for $y_1=1$.

Pair $(S,T)$ of spheres is called a {\it strong realization} of pair $(\vec x, \vec y)$.

{\it The `if' part in Theorem 1} is implied by the following Theorem.

\smallskip
{\bf Theorem 1'.}
{\it Each pair of tree-like sequences is strongly realizable.}
\smallskip

In order to prove this Theorem we need the following Claim.

\smallskip
{\bf Lemma 1.}
{\it Let $\vec x = (x_1, x_2, \dots, x_n), \vec y = (y_1, y_2, \dots, y_n)$ be tree-like sequences in which all the units are situated at the end.
Assume that $x_1\ge y_1$ and denote 
$$\vec x':=(x_1-y_1+1, x_2, x_3, \dots, x_{n-y_1+1})\quad\text{and}\quad\vec y':=(y_2, y_3, \dots, y_{n-y_1+2}).$$

(a) Then sequences $\vec x'$ and $\vec y'$ are tree-like.

(b) If pair of sequences $(\vec x',\vec y')$ is strongly realizable, then pair $(\vec x,\vec y)$ is strongly realizable.}

 \comment

$\bullet$ $\vec {x'}=(a, x_1, x_2, \dots, x_n, 1,1, \dots,1)$ (the number of new $1$'s is $a-2$),
\linebreak
$\vec {y'}=(y_1+a-1, y_2, y_3, \dots, y_n, 1, 1, \dots, 1)$ (the number of new $1$'s is $a-1$).

$\bullet$ $\vec {x'}=(a, x_1, x_2, \dots, x_n,1,1, \dots,1)$ (the number of new $1$'s is $a-2$),
\linebreak
$\vec {y'}=(a, y_1, y_2, y_3, \dots, y_n,1,1, \dots,1)$ (the number of new $1$'s is $a-2$).

Then each pair of tree-like sequences is realizable.

In fact in the proof of Lemma 1 we construct the  pair $(\{x\}, \{y\})$ using the pair $(\{x'\}, \{y'\})$.


\endcomment

\smallskip
{\it Proof.}
(a)
Let $s$ be a number of units in $\vec x$.
Then
$$2n-2=x_1 +\dots +x_n\ge x_1+2(n-1-s)+s=2n-2+x_1-s.$$
So $s\ge x_1$.
Analogously $r\ge y_1$, where $r$ is a number of units in $\vec y$.

Then
$$x_{n-y_1+1}=x_{n-y_1+2}=\dots=x_{n-y_1+1}=\dots=x_n=y_{n-y_1+1}=y_{n-y_1+2}=\dots=y_n=1.$$
$$\text{Hence}\quad(\sum_{i=1}^{n-y_1+1}x_i)-y_1+1=(\sum_{i=1}^nx_i) - y_1+1-(y_1-1)=2(n-y_1+1)-2$$
$$\text{and}\quad(\sum_{i=2}^{n-y_1+2}y_i)=(\sum_{i=1}^ny_i) - y_1-(y_1-2)=2(n-y_1+1)-2.$$
So sequences $\vec x'$ and $\vec y'$ are  tree-like. QED

\smallskip
(b)
 Take spheres $S', T'$ realizing pair $(\vec x',\vec y')$  of sequences.
Take a circle of $S' \cap T'$ from condition (2).
This circle bounds

$\bullet$ in $S'-T'$ a connected component, say $C$, that has $x_1-y_1+1$ neighbors,

$\bullet$  in $T'-S'$ a disk, say $D$.

\begin{figure}[h]
  \centerline{\includegraphics[width=10cm]{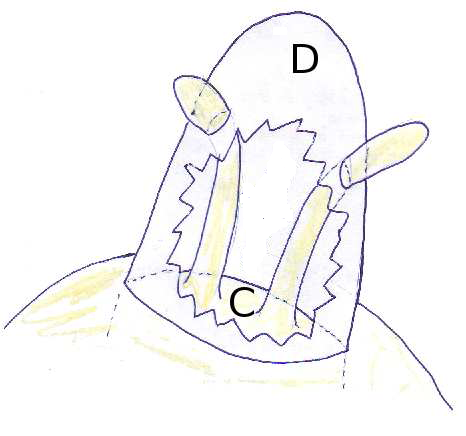}}
\caption{Inductive construction}\label{fingers}
\end{figure}

We modify  spheres $S', T'$ by joining
$C$ and $D$  by $y_1-1$ fingers, see Figure \ref{fingers}.
 Denote the new  spheres by $S$ and $T$.
Let  us prove that they realize pair $(\vec x, \vec y)$ of sequences.

Condition (1) is satisfied for $S,T$ because

$\bullet$ each component of $S'-T'$ except $C$
is also a component of $S-T$,

$\bullet$ $C$ is separated by $y_1-1$ circles of $(S\cap T) - (S'\cap T')$  into $y_1-1$ disks and a component with $(x_1-y_1+1)+(y_1-1)=x_1$ neighbors.

and

$\bullet$ each component of $T'-S'$ except $D$ is also a component of $T-S$,

$\bullet$ $D$ is separated by $y_1-1$ circles into $y_1-1$ disks and a component with $y_1$ neighbors.

Any circle of $(S\cap T) - (S'\cap T')$ satisfies condition (2).
QED

\smallskip
{\it Proof of Theorem 1'.}
By induction on the length $n$ of the sequences.
For each tree-like sequence of $n$ numbers we have $n\ge2$.
The induction base is $n=2$ and is clear.

Let us prove the induction step.
Suppose Theorem 1' is proved for $2,3, \dots, n-1\ge 2$.
 Let us prove it for $n$.

We can reorder our sequences so that 

$\bullet$ the $1$'s will be at the ends, 

$\bullet$ if $x_1>1$ then $x_1$ will not change its position, and 

$\bullet$ if $y_1>1$ then $y_1$ will not change its position. 

By Lemma 1 (a), (b) and by the induction hypothesis the `reordered' pair is strongly realizable.

Let us prove that the initial pair is strongly realizable. 
Take spheres $S,T$ realizing the new sequences.
So $S,T$ satisfy condition (1) from the definition of the strong realizability for the initial sequences.
Also,

$\bullet$ if $x_1>1$, then conditions (2) for the reordered and for the initial sequences $\vec x$ are equivalent;

$\bullet$ if $x_1=1$, then the circle from condition (2) for the reordered sequence $\vec x$ bounds a disk, so the circle bounds a component with $x_1=1$ neighbor.

Same holds for $\vec x$ replaced by $\vec y$.
So condition (2) is also satisfied for the initial sequences.
Thus $S,T$ strongly realize the initial sequences.
QED

\smallskip

In order to prove Theorem 2 we need the following elementary lemma.

\smallskip
{\bf Lemma 2.} {\it Let $x_1\ge  x_2\ge  \dots\ge x_n$ be a tree-like sequence. Let $p$, $q$ be positive integers such that $p \ge q>1$ and $p+q=n+1$.
Then there exist two tree-like sequences $a_1, a_2, \dots, a_p$ and $b_1, b_2, \dots , b_q$ such that $a_1+b_1=x_1$ and ordered sets $(a_2, a_3, \dots, a_p, b_2, b_3, \dots, b_q)$ and $(x_2, x_3, \dots, x_n)$ are the same up to reordering.}

\smallskip

{\it Proof.}
Let $r=r(\vec x)$ be the number of those $x_i$'s that are greater than 1.
Let $z_s=x_2+x_3+ \dots +x_s$.
For each $s\le r$ let
$$a_1=p-(z_s -s+3)+1,\quad a_i=x_i\quad\text{for}\quad 2\le i\le s\quad\text{and}\quad a_i=1\quad\text{for}\quad s+1\le i\le p,$$
$$b_1=x_1-a_1,\quad b_i=x_{i+s-1}\quad\text{for}\quad 2\le i\le r-s+1,\quad b_i=1\quad\text{for}\quad r-s+2\le i\le q=n+1-p.$$
Since $s\le r$, the sequence $b_1, b_2, \dots , b_q$ is well-defined.
For each $i$ we have that $a_i$ and $b_i$ depend on $s$.

We have
$$a_1+a_2+ \dots +a_p = p-(z_s -s+3)+1+z_s+p-s=2p-2,$$
 i.e. the sequence $a_1, a_2, \dots , a_p$ is tree-like.
Also
$$b_1+b_2+ \dots +b_q = z_n-a_1-a_2- \dots -a_p=2n-2-2p+2=2q-2,$$
 i.e. the sequence $b_1, b_2, \dots, b_q$ is tree-like.

It remains to prove that there exists $s\le r$ such that $1\le a_1 \le x_1-1$.
For each $i<r$ we have $x_1 \ge x_i$, so
$$z_i-i+x_1+1 \ge (z_{i+1}-(i+1)+3)-1.$$
In other words,

\centerline {$2=z_1-1+3$,}

\centerline {$z_1-1+x_1+1 \ge (z_2-2+3)-1$, }

\centerline {$z_2-2+x_1+1 \ge (z_3-3+3)-1$, }

\centerline {\dots, }

\centerline {$z_{r-1}-(r-1)+x_1+1 \ge (z_r-r+3)-1$,}

\centerline {$z_r-r+x_1+1=n-1$.}

Here the last equality is not
analogous to the previous equalities but follows because sequence $x_1, x_2, \dots, x_n$ is tree-like and $1=x_{r+1}=\dots=x_n$.
Since $2\le p\le n-1$, there exists $s\le r$ such that
$$z_s -s+3\le p\le z_s-s+x_1 +1 \quad\Leftrightarrow\quad 1\le a_1 \le x_1-1.\quad QED$$
\smallskip

{\it Proof of the `if' part in Theorem 2.}
Let $$m_1:=(n_2+n_3-n_1+1)/2,\quad m_2:=(n_1+n_3-n_2+1)/2,\quad m_3:=(n_1+n_2-n_3+1)/2.$$
So $$m_1+m_2=n_3+1,\quad m_1+m_3=n_2+1,\quad m_2+m_3=n_1+1.$$
Hence by Lemma 2 there exist sequences
$$
p_{1 1} , p_{1 2} , \dots , p_{1 m_3},\quad p_{2 1} , p_{2 2} , \dots , p_{2 m_1},
\quad p_{3 1} , p_{3 2} , \dots , p_{3 m_2},
$$
$$
q_{1 1} , q_{1 2} , \dots , q_{1 m_2},\quad q_{2 1} , q_{2 3} , \dots , q_{2 m_3},
\quad q_{3 1} , q_{3 2} , \dots , q_{3 m_1},
$$
such that $p_{k-1, 1}+q_{k+1, 1}=x_{k 1}$ and
ordered sets $$(p_{k-1, 2}, p_{k-1, 3}, \dots, p_{k-1, m_{k+1}}, q_{k+1, 2}, q_{k+1, 3}, \dots, q_{k+1, m_{k-1}})\quad \text{and}\quad (x_{k 2}, x_{k 3}, \dots, x_{k n_k})$$ are the same up to reordering.


By Theorem 1' there exist  spheres
$$Q_1, P_1, Q_2, P_2, Q_3, P_3\subset\R^3\quad\text{such that}\quad Q_k \cap Q_{k+1}=\varnothing,\quad Q_k\cap P_l=\varnothing\quad\text{if}\quad l\neq k-1\quad\text{and}$$

$\bullet$ $Q_k-P_{k-1}$ is a disjoint union of $m_{k+1}$ connected components, $i$-th one has $q_{k i}$ neighbors;


$\bullet$ $P_{k-1}-Q_k$ is a disjoint union of $m_{k+1}$ connected components, $i$-th one has $p_{k-1, i}$ neighbors;


$\bullet$ the boundary of some connected component of $\R^3-P_{k-1}-Q_k$ contains a component, say $\widetilde q_k$, with $q_{k, 1}$ neighbors on $Q_k$ and a component, say $\widetilde p_{k-1}$, with $p_{k-1, 1}$ neighbors on $P_{k-1}$.

\begin{figure}[h]
  \centerline{\includegraphics[width=10cm]{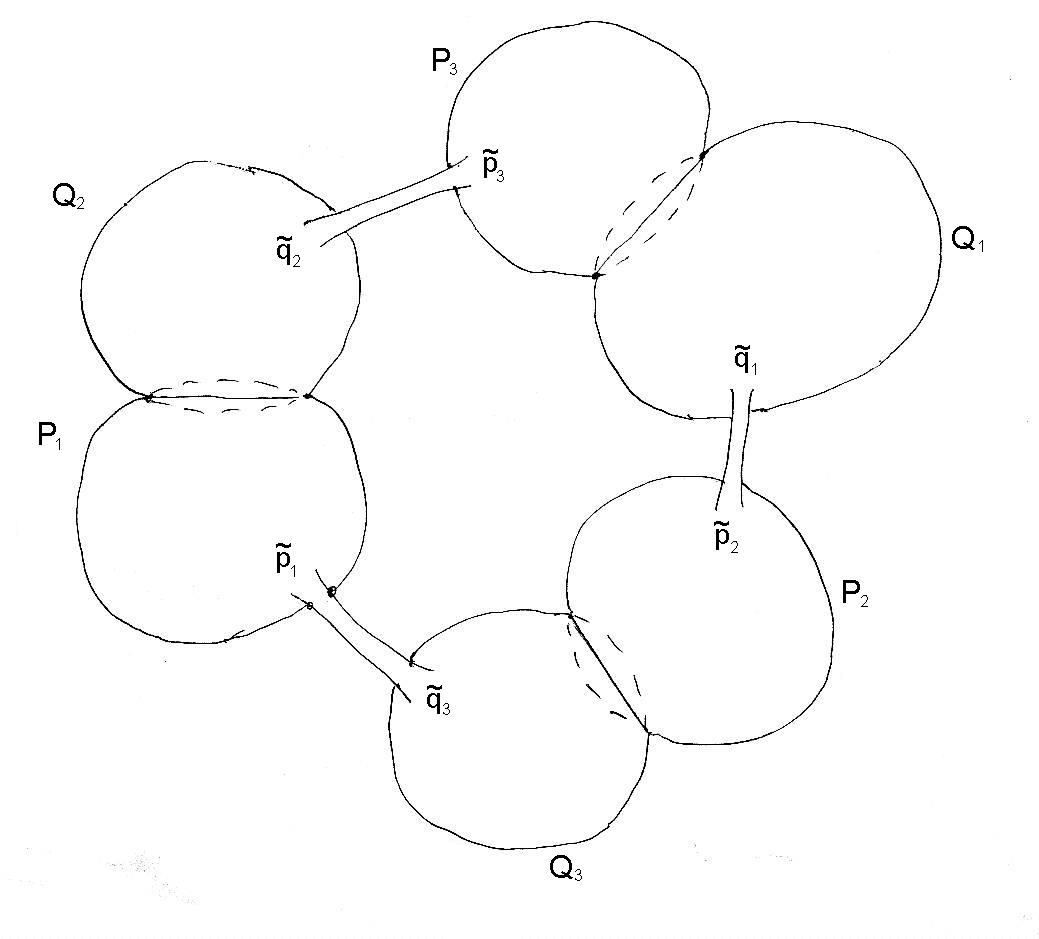}}
\caption{Construction for Theorem 2}\label{3spheres}
\end{figure}
We may assume that in Theorem 1' the boundary of the infinite component of $\R^3-S-T$ contains a circle from the condition 2, because applying an inversion with respect to a point in the component whose boundary contains this circle we obtain required spheres. The pairs of spheres $P_{k-1}, Q_k$ are assumed to be such.

For $1\le k \le 3$ let $f_k$ be the connected sum of  spheres $Q_{k+1}$ and $P_{k-1}$ along a small tube joining the two components $\widetilde q_{k+1}$ and $\widetilde p_{k-1}$, see Figure \ref{3spheres}.
This can be done without intersections of the three tubes.

Then $f_k-f_{k+1}-f_{k+2}$ is as required for each $1\le k\le 3$. QED.


\bigskip
{\bf 3. Some alternative proofs and improvements.}

\smallskip
{\it Scheme of an alternative proof of Theorem 1'.}
One can prove that

$\bullet$ pair $(\vec x,\vec x)$ is strongly realizable for each tree-like sequence $\vec x$;

$\bullet$ if pair $(\vec x,\vec y)$ is strongly realizable, then the following pair $(\vec {x'}, \vec {y'})$ is strongly realizable:

$\vec {x'}=(a, x_1, x_2, \dots, x_n, 1, 1, \dots, 1)$ (the number of new $1$'s is $a-2$),

$\vec {y'}=(y_1+a-1, y_2, y_3, \dots, y_n, 1, 1, \dots, 1)$ (the number of new $1$'s is $a-2$).

 Then each pair of tree-like sequences is strongly realizable.

\smallskip
{\it Proof of the strong realizability of pair $(\vec x, \vec x)$.} 
Let $f$ be the unite cube.
Take a family of circles on $f$ `realizing' $\vec x$.
(The existence of such a family is proved by induction involving deleting the unity entry from $\vec x$.)
Color the complements in $f$ to these circles into black and white so that neighboring components have different colors.
Take a sphere $g$ close to $f$ and such that

$\bullet$ $f\cap g$ is the disjoint union of the circles;

$\bullet$ each black component of $g$ is inside $f$;

$\bullet$ each white component of $g$ is outside $f$.
QED

\smallskip
{\it Yet another alternative proof of Theorem 1'} could possibly be obtained by assuming that a pair $(\vec x, \vec y)$ is realizable and realizing the pair $(\vec{x'},\vec{y'})$ for
$$\vec{x'}\in\{(1,x_1+1,x_2,x_3,\dots,x_n),(x_1+1,x_2,x_3,\dots,x_n,1)\}\quad\text{and}$$
$$\vec{y'}\in\{(1,y_1+1,y_2,y_3,\dots,y_n),(y_1+1,y_2,y_3,\dots,y_n,1)\}.$$

\smallskip
{\it Remark: an improvement of Theorem 1.}
Define of a sketch $A(\vec x)$ for a tree-like sequence $\vec x=(x_1, x_2,\dots , x_n)$, where $x_i=1$ if and only if $i>r(x_1,x_2,\dots,x_n)$.
Draw $r+1$ circles on sphere $S^2$ so that
these circles split $S^2$ into two disks and $r$ annuli (an annulus is a disk with one hole).
We call them {\it main circles}.
For each $i$ from 1 to $r$ draw $x_i-2$  non-intersecting disks in the $i$-th annulus from the top.
We have drawn $x_1+x_2+\dots+x_{r}-r+1=n-1$ circles.

Recall that the sum of zero summands is zero, so that $x_1+x_2+\dots+x_i=0$ for $i=0$.
Denote all these $n-1$ circles as follows:


$\bullet$ for each $i$ from 0 to $r$, the $(i+1)$-st from the top main circle we denote
$A_{x_1+x_2+\dots+x_i-i+1}$;


$\bullet$ for each $i$ from 0 to $r-1$ circles in the interior of the $(i+1)$-st annulus from the top we denote
$$A_{x_1+x_2+\dots+x_i-i+2},\quad A_{x_1+x_2+...+x_i-i+3},\quad  \dots,\quad A_{x_1+x_2+\dots+x_i+x_{i+1}-(i+1)};$$


The ordered set $A(\vec x)=(A_1, A_2, \dots,A_{n-1})$ of circles on $S^2$ is called {\it the sketch for sequence $\vec x=(x_1, x_2,\dots, x_n)$}.

\smallskip
{\bf Improvement of Theorem 1. } {\it Let $A_1, A_2, \dots,A_{n-1}$ be the sketch for a tree-like sequence $\vec x = (x_1, x_2,\dots, x_n)$, in which all the units are at the end. Also let $B_1, B_2, \dots,B_{n-1}$ be the sketch for a tree-like sequence $\vec y = (y_1,  y_2, \dots,  y_n)$, in which all the units are at the end. Then there exist two embeddings $F, G:S^2\to \R^3$ such that
$$F(A_i)=G(B_i)\quad\mbox{for each}\quad 1\le i\le n-1 \quad\mbox{and}\quad F(S^2)\cap G(S^2)=\sqcup_{i=1}^{n-1} F(A_i).$$
}

\comment
{\it Proof of the sufficiency in Theorem 1.}
We may assume that for given tree-like sequences $x_i=1$ for $i>r(x_1,x_2,\dots,x_n)$ and $y_i=1$ for $i>r(y_1,y_2,\dots,y_n)$.
Take sketches $A_1, A_2, \dots, A_{n-1}$ and $B_1, B_2, \dots,B_{n-1}$ corresponding to the sequences.
Then $S^2-A_1-A_2-\dots-A_{n-1}$ is the disjoint union of $n$ connected components,  which can be numbered so that the $i$-th component has $x_i$ neighbors.
Analogous property holds for the sequence of $y$'s and the corresponding sketch.
Now the sufficiency in Theorem 1 is implied by the following result.
QED
\endcomment

\smallskip

{\it Idea of an alternative proof of Lemma 2.}
A referee suggested the following interpretation of our proof of Lemma 2 in terms of sketch. Take an integer $s$ (see details in the formal proof). Realize $x_2, x_3, \dots, x_s, x_1, x_{s+1}, \dots, x_n$ by a sketch. Between $s$-th and $s+1$-th main circles draw an additional circle that split the sketch into two subsketches of $p-1$ and $q-1$ circles satisfying to the conditions of Lemma 2.

\smallskip
{\bf Acknowledgments.}
The author is grateful for useful remarks and discussions to professor S. Lando and to an anonymous referee of Moscow Mathematical Conference of High-School Students.

\bigskip
{\bf References}

[A] S. Avvakumov, A counterexample to the Lando conjecture on
intersection of spheres in 3-space, preprint, to be put in arxiv, 2012.

[H] T. Hirasa, Dissecting the torus by immersions,
Geometriae Dedicata, 145:1 (2010), 33-41.

[N] T. Nowik, Dissecting the 2-sphere by immersions,
Geometriae Dedicata 127, (2007), 37-41.
http://arxiv.org/abs/math/0612796

[T] S. Avvakumov, A. Berdnikov, A. Skopenkov, A. Rukhovich, How do curved spheres intersect in 3-space, or 2-dimensional meandra, http://www.turgor.ru/lktg/2012/3/3-1en\_si.pdf

[Tr] S. Avvakumov, A. Berdnikov, A. Skopenkov, A. Rukhovich, How do curved spheres intersect in 3-space, or 2-dimensional meandra, in russian, http://www.turgor.ru/lktg/2012/3/3-1ru\_si.pdf

\end{document}



\comment


 \normalsize

\endabstract

\endcomment

{\centerline{\bf Как пересекаются многогранники в трехмерном пространстве?}}
\bigskip
{\centerline{\bf Рухович А.}}
\bigskip
{\bf Введение.}

Как может быть устроено пересечение многогранников в трехмерном пространстве?
Приведем формализацию этой проблемы, восходящую к Т. Новику 2007 [N, H].
Определение многогранника (не обязательно выпуклого) приведено в [D], см. также [W].

\begin{figure}[h]
  \centerline{ \includegraphics[width=6cm]{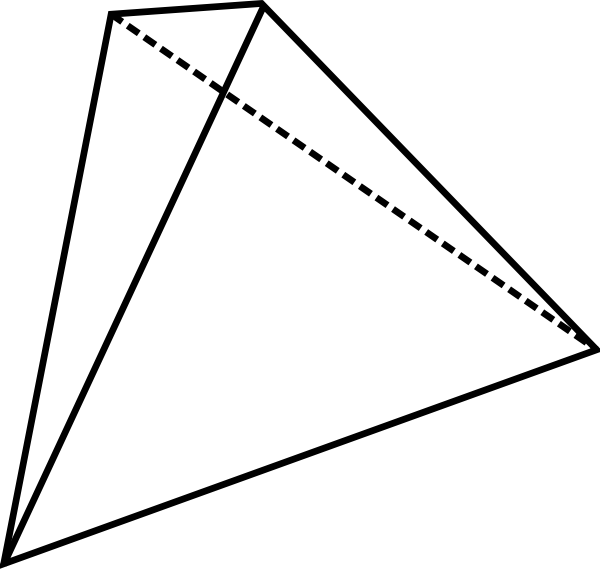} \includegraphics[width=6cm]{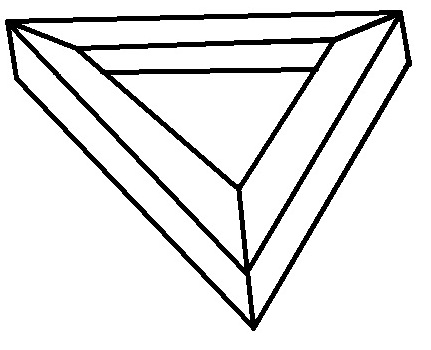}}
 {  Многогранник, топологически эквивалентный сфере (слева),
  и не топологически эквивалентный сфере (справа)}
\end{figure}

Будем рассматривать многогранники, топологически эквивалентные сфере (т.е. гомеоморфныe сфере) [D], см. также [W].

\begin{figure}[h]
  \centerline{ \includegraphics[width=10cm]{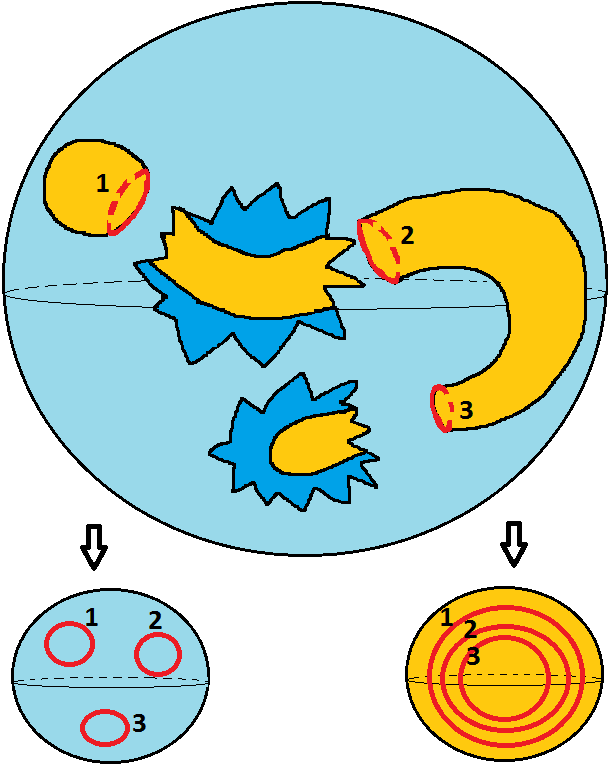}}
\end{figure}

Для упрощения рисунков будем изображать вместо многогранника `близкую' к нему {\it криволинейную} поверхность.
Например, сферу или `сосиску', как на вышеприведенном рисунке.
На этом
рисунке два многогранника $f$ и $g$ пересекаются по трем окружностям.
Эти окружности разбивают каждый многогранник на четыре части.
На первом многограннике у первой части 3 соседей, у каждой из остальных частей --- по 1 соседу.
На втором многограннике у первой и четвертой части по 2 соседа, у второй и третьей части --- по 1 соседу.

Объединение частей на первом многограннике является дополнением $f-g$ (другое обозначение: $f\backslash g$); аналогичный смысл имеет обозначение $g-f$.


\bigskip
{\bf Формулировка основных результатов.}

{\it Парой многогранников общего положения} называется пара несамопересекающихся многогранников в трехмерном пространстве,
топологически эквивалентных сфере,
пересечение которых состоит из конечного числа замкнутых несамопересекающихся ломаных.

Дадим формальное определение соседства (надеемся, что следующий текст будет понятен и без него).
Две связные компоненты дополнения $f-g$ назовем {\it соседями в $f$}, если их замыкания пересекаются.
Аналогично определим {\it соседей в $g$} для компонент дополнения $g-f$.




\smallskip
{\bf Теорема 1.} {\it Пусть $n$ натуральное число и $x_1,x_2,\dots,x_n$, $ y_1,y_2, \dots, y_n$ последовательности натуральных чисел.
Тогда существует пара многогранников $f$ и $g$ общего положения таких, что

$\bullet$
$f-g$ имеет $n$ компонент связности, $i$-тая из которых имеет $x_i$ соседей из $f$;

$\bullet$
$g-f$ имеет $n$ компонент связности $i$-тая из которых имеет $y_i$ соседей из $g$;

тогда и только тогда, когда}
$$\sum_{i=1}^nx_i=\sum_{i=1}^ny_i=2n-2.$$

 Часть `только тогда' известна и доказана в [N]
(здесь представлены два простых доказательства, одно из которых использует теорему Жордана, а второе --- доказательство T. Новика).
Часть `тогда'
новая.
Она особенно интересна в связи с отрицательным решением следующей гипотезы.

\smallskip
{\bf Гипотеза Ландо.} {\it Пусть $A$ объединение непересекающихся окружностей на $S^2$.
Аналогично пусть $B$ объединение такого же числа непересекающихся окружностей на $S^2$.
Тогда существуют многогранники $f,g\subset\R^3$, топологически эквивалентные сфере, такие что $f\cap g$ объединение ломаных, аналогичное $A$ на $f$ и $B$ на $g$.}

\begin{figure}[h]
  \centerline{\includegraphics[width=12cm]{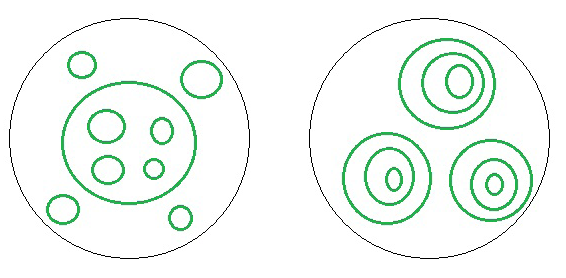}}
  \centerline{  Контрпример к гипотезе Ландо.}
\end{figure}

Контрпример (см. рисунок выше) найден и доказан Сергеем Аввакумовым [A].

Поясним, как теорема 1 связана с классической тематикой изучения {\it степенных последовательностей} графов.
Построим граф $F$, вершины которого --- связные компоненты дополнения $f-g$;
ребрами соединим вершины, если соответствующие компоненты --- соседи.
Аналогично определим граф $G$.
Тогда Теорема 1 описывает пары степенных последовательностей таких графов.

\begin{figure}[h]
  \centerline{\includegraphics[width=15cm]{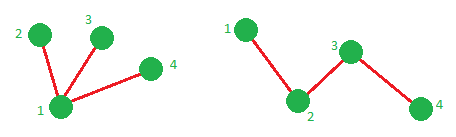}}
  \centerline{  Графы $F$ и $G$ для многогранников с рисунка `пример'}
\end{figure}

\smallskip
Сформулируем аналог теоремы 1 для пересечения {\it трех} многогранников.
Назовем последовательность натуральных чисел $x_1, x_2, \dots, x_n$ {\it деревянной} если $\sum_{i=1}^nx_i=2n-2$.

\smallskip
{\bf Теорема 2.} {\it Пусть $n_1, n_2, n_3$ натуральные числа и
  $$x_{1 1}, x_{1 2}, \dots, x_{1 n_1}, x_{2 1}, x_{2 2}, \dots ,x_{2 n_2}, x_{3 1}, x_{3 2}, \dots, x_{3 n_3}$$
последовательности натуральных чисел. Тогда существуют многогранники $f_1, f_2, f_3\subset\R^3$,
любые два из которых образуют пару многогранников общего положения, причем

$\bullet$ $f_1 \cap f_2 \cap f_3 =\varnothing$;

$\bullet$ $f_k-f_{k+1}-f_{k+2}$ имеет $n_k$ компонент связности, $i$-тая из которых имеет $x_{k i}$ соседей в $f_k$, для любого $1\le k \le 3$

тогда и только тогда, когда последовательности деревянные, $n_1+n_2+n_3$ нечетно и $n_k < n_{k+1}+n_{k+2}$ для любого $1\le k \le 3$.}

\smallskip
Здесь индексы $k, k+1, k+2$ рассматриваются по модулю 3.

Часть `только тогда' тривиальна. Часть `тогда' 
новая.

\bigskip
{\bf Гипотезы.}

Каков аналог теорем 1 и 2 для пересечения {\it более, чем трех} многогранников?

\smallskip
{\bf Гипотеза 1.}{\it Заданы натуральные числа $n_1, n_2, n_3, \dots, n_k$ и $k$ последовательностей
$$x_{1 1}, x_{1 2}, \dots, x_{1 n_1},\quad x_{2 1}, x_{2 2}, \dots ,x_{2 n_2},\quad \dots, \quad x_{k 1}, x_{k 2}, \dots, x_{k n_k}$$
натуральных чисел. Тогда существуют $k$ многогранников  $f_1, f_2, \dots,  f_k\subset\R^3$, топологически эквивалентных сфере, таких, что

$\bullet$ никакие три из них не пересекаются в одной точке;

$\bullet$  для каждого $1\le i \le k$ и $1\le j \le n_i$ выполнено $f_i-f_{i+1}-f_{i+2}-\dots - f{i+k-1}$ имеет $n_i$ компонент связности, $j$-тая из которых имеет $x_{i j}$ соседей в $f_i$

тогда и только тогда когда

$\bullet$ все $k$ последовательностей деревянные,

$\bullet$ $(n_1-1)+(n_2-1)+\dots+(n_k-1)$ четно и

$\bullet$ $2(n_i-1) \le (n_1-1)+(n_2-1)+\dots+(n_k-1)$ для каждого $1\le i \le k$.}

\smallskip

Каков аналог теоремы 1 для пересечения многогранников {\it  топологически эквивалентных сферам с ручками}?

\smallskip
{\bf Гипотеза 2.}{\it Заданы натуральные числа $g_1, g_2, n_1, n_2$ и две последовательности $$x_{1 1}, x_{1 2}, \dots, x_{1 n_1},\quad x_{2 1}, x_{2 2}, \dots, x_{2 n_2}$$ натуральных чисел. Тогда существуют многогранники  $f_1, f_2\subset\R^3$ таких что

$\bullet$ $f_1$ топологически эквивалентен сфере с $g_1$ ручками,  $f_2$ топологически эквивалентен сфере с $g_2$ ручками;


$\bullet$ $f_i-f_{i+1}$ имеет $n_i$ компонент связности, $j$-тая из которых имеет $x_{i j}$ соседей в $f_i$, для каждого $1\le i \le k$

тогда и только тогда когда  $s = x_{1 1} + x_{1 2} + \dots + x_{1 n_1} = x_{2 1} + x_{2 2} + \dots + x_{2 n_2}$ - четное число, $2\ge s - 2n_k \ge  2-2g_k$.}

\smallskip
{\bf Проблема.} {\it Какими могут быть наборы количеств пересечений у компонент, если разрешить пересечения трех 
многогранников в одной точке?}

\smallskip
Было бы интересно решить аналогичную проблему при наличии самопересечений сфер.
Интересны оба варианта --- с точками троекратного самопересечения или без.

\comment
{\bf Замечание.} Следующий занумерованный аналог Гипотезы Ландо неверен.

 {\it
 Пусть $A_1, A_2, \dots, A_n$ непересекающиеся окружности на $S^2$.
Аналогично пусть $B_1, B_2, \dots, B_n$ непересекающиеся окружности на $S^2$.
Тогда существуют вложения $f,g: S^2 \to \R^3$ такие что
$f(B_s)=g(A_s)$ для любого $s \in \{1, 2, \dots,n \}$ и $f(S^2), g(S^2)$ не имеют других точек пересечения.}

 {\it Контрпример для $n=3$.}
Пусть $A_1$ окружность на один радиан северной широты, $A_2$ экватор сферы (т.е., ноль радиан северной широты),
$A_3$  окружность на один радиан южной широты.
Пусть $B_1=A_2$, $B_2=A_1$, $B_3=A_3$.
\endcomment

%


\bigskip
{\bf Доказательства.}


{\it Доказательство части `только тогда ' в Теореме 1.}
Напомним определение графа $F$.
Вершины $F$ это компоненты связности $f-g$.
Две вершины соединим ребром если соответствующие компоненты - соседи. Тогда количество вершин равно $n$. Количество ребер равно количеству окружностей в $f\cap g$. Это число равно 
$\sum_{i=1}^nx_i/2$. Очевидно что $F$ связен. По теореме Жордана $F$ разбивается любой вершиной на различные компоненты связности. Значит $F$ дерево. Тогда количество его ребер равно $n-1=\sum_{i=1}^nx_i/2$. QED

 \smallskip
{\it Доказательство части `только тогда' Теоремы 1 представленное T. Nowik`ом.}
Докажем индукцией по количеству окружностей. Утверждение верно для одной окружности (всего 2 диска на каждой сфере и $n=2$). Каждая новая окружность делит одну компоненту на две и добавляет две граничные окружности.
QED

{\it Доказательство части `только тогда' в Теореме 2.}
Необходимость первого условия очевидно следует из части `только тогда' Теоремы 1.

Пусть $m_3, m_2, m_1$ количества окружностей в $f_1\cap f_2$, $f_1\cap f_3$ and $f_2\cap f_3$.
Тогда $n_1=m_3+m_2+1$, $n_2=m_3+m_1+1$, $n_3=m_2+m_1+1$.

Значит $n_1+n_2+n_3=2(m_3+m_2+m_1)+3$ нечетно.

Поскольку $2m_k+1>0$ получаем $n_k < n_{k+1}+n_{k+2}$ для любого $1\le k\le 3$.
QED



Для доказательства Теоремы 1 введем следующие определения.

\smallskip
Пара $\vec x = (x_1,x_2,\dots ,x_n), \vec y = (y_1,y_2,\dots ,y_n)$ последовательностей положительных целых чисел называется {\itshape реализуемой}, если существуют две криволинейные сферы $S$ и $T$ в трехмерном пространстве, такие что

(1) их пересечение состоит из $n-1$ окружностей и делит

\quad  $\bullet$ $S$ на $n$ связных компонент, которые можно занумеровать так, что у $i$-ой связной компоненты $x_i$ соседей в $S$, играницей

\quad  $\bullet$ $T$ на $n$ связных компонент, которые можно занумеровать так, что у $i$-ой связной компоненты $y_i$ соседей в $T$;

(2) в $S\cap T$ существует окружность, являющаяся границей некоторых диска и компоненты с $x_1$ соседями в $S-T$,
а также границей некоторых диска и компоненты с $y_1$ соседями в $T-S$.

Пара $(S,T)$ сфер называется {\itshape реализацией} пары $(\vec x, \vec y)$.

\smallskip
Часть `тогда' в Теореме 1 вытекает из следующего 
факта.

\smallskip
{\bf Теорема 1'.} {\it Любая пара деревянных последовательностей реализуема.}

Эта теорема выводится из следующей Леммы.
\smallskip

{\bf Лемма 1.} {\it Пусть $\vec x, \vec y$ --- деревянные последовательности, в которых все единицы находятся в конце, и $x_1\geq y_1$. Тогда

(a) Последовательности $\vec x' := (x_1-y_1+1,x_2,x_3,\dots ,x_{n-y_1+1})$ и
$\vec y' = (y_2,y_3,\dots ,y_{n-y_1+2})$ --- деревянные.

(b) Если пара последовательностей $\vec x':=(x_1-y_1+1, x_2, x_3, \dots, x_{n-y_1+1})$,
$\vec y':=(y_2, y_3, \dots, y_{n-y_1+2})$ реализуема, то и пара $(\vec x,\vec y)$ тоже.
}

{\it Доказательство Леммы 1.}
(а)
Пусть $s$ --- количество единиц в последовательности $x_1, x_2, \dots, x_n$.
Тогда $2n-2=x_1 +\dots +x_n\ge x_1+2(n-1-s)+s=2n-2+x_1-s$.
Значит, $s\ge x_1$.
Тогда
$$x_{n-y_1+1}=x_{n-y_1+2}=\dots=x_{n-y_1+1}=\dots=x_n=y_{n-y_1+1}=y_{n-y_1+2}=\dots=y_n=1.$$
$$\text{Следовательно}\quad(\sum_{i=1}^{n-y_1+1}x_i)-y_1+1=(\sum_{i=1}^nx_i) - y_1+1-(y_1-1)=2(n-y_1+1)-2$$
$$\text{и}\quad(\sum_{i=2}^{n-y_1+2}y_i)=(\sum_{i=1}^ny_i) - y_1-(y_1-2)=2(n-y_1+1)-2.$$
Значит новые последовательности деревянные.

(b)
Рассмотрим сферы $S', T'$, реализующие пару $(\vec x',\vec y')$ последовательностей.
Рассмотрим окружность в пересечении $T' \cap S'$, из условия (2) в определении сильной реализуемости.
Это окружность, ограничивающая

$\bullet$ в $S'-T'$ связную компоненту, назовем ее $C$, у которой $x_1-y_1+1$ соседей,

$\bullet$ в $T'-S'$ --- диск, назовем его $D$.

Изменим сферы $S', T'$ соединением компонент
$C$ и $D$ при помощи $y_1-1$ `пальцев', см. рис..
 Новые сферы обозначим через $S, T$.
Докажем, что они реализуют пару  $(\vec x, \vec y)$ последовательностей.

Условие (1) выполнено для $S, T$, поскольку

$\bullet$ каждая компонента $S'-T'$, кроме $C$,
также является и компонентой дополнения $S-T$,

$\bullet$ $C$ делится $y_1-1$ окружностями из $(S\cap T) - (S'\cap T')$ на $y_1-1$ дисков и компоненту с $(x_1-y_1+1)+(y_1-1)=x_1$ соседями,

и

$\bullet$ каждая компонента дополнения  $T'-S'$, кроме $D$, также является и компонентой дополнения  $T-S$,

$\bullet$ $D$ разделено $y_1-1$ окружностями на $y_1-1$ дисков и компоненту с $y_1$ соседями.

Любая окружность из $(S\cap T) - (S'\cap T')$ удовлетворяет условию (2).

\begin{figure}[h]
  \centerline{\includegraphics[width=15cm]{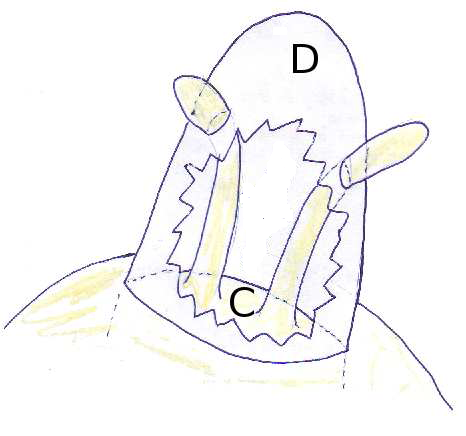}}
  \centerline{ Индуктивное построение.}
\end{figure}


Из Леммы 1 также следует обобщение Теоремы 1, нужное для доказательства Теоремы 2.

\smallskip
{\bf Теорема 1'.} {\it Пусть $n$ натуральное число и $x_1, x_2, \dots, x_n$, $y_1, y_2,  \dots, y_n$ деревянные последовательности натуральных чисел.
Тогда существует пара многогранников $f,g$ общего положения такая, что

$\bullet$
$f-g$ имеет $n$ компонент связности, $i$-тая из которых имеет $x_i$ соседей из $f$;

$\bullet$
$g-f$ имеет $n$ компонент связности, $i$-тая из которых имеет $y_i$ соседей из $g$;

$\bullet$ граница некоторой связной компоненты $\R^3-f-g$ содержит компоненту с $x_1$ соседями на $f$ и компоненту с $y_1$ соседями на $g$.

}



Для доказательства части `тогда' в теореме 2 используем следующую элементарную лемму.

\smallskip
{\bf Лемма 2.} {\it Пусть $x_1\ge  x_2\ge  \dots\ge x_n$ деревянная последовательность. Пусть $p$, $q$ натуральные числа, что $p \ge q>1$ и $p+q=n+1$.
Тогда существуют деревянные последовательности $a_1, a_2, \dots, a_p$ и $b_1, b_2, \dots , b_q$ такие что $a_1+b_1=x_1$ и упорядоченные наборы  $(a_2, a_3, \dots, a_p, b_2, b_3, \dots, b_q)$ и $(x_2, x_3, \dots, x_n)$ одинаковы с точностью до перестановки.}

\smallskip

{\it Доказательство.}
Пусть $r=r(x_1, x_2,\dots, x_n)$.
Пусть $z_s=x_2+x_3+ \dots +x_s$.
Для каждого $s\le r$ пусть
$a_1=p-(z_s -s+3)+1, a_i=x_i$ для всех $2\le i\le s$ и $a_i=1$ для всех $s+1\le i\le p,$
$b_1=x_1-a_1, b_i=x_{i+s-1}$ для всех $2\le i\le r-s+1, b_i=1$ для всех $r-s+2\le i\le q=n+1-p.$
Поскольку $s\le r$, последовательность $b_1, b_2, \dots , b_q$ определена непротиворечиво.
Для всех $i$ имеем что $a_i$ и $b_i$ зависят от $s$.

Получаем
  $a_1+a_2+ \dots +a_p = p-(z_s -s+3)+1+z_s+p-s=2p-2,$
  т.е. последовательность $a_1, a_2, \dots , a_p$ деревянная.
Также
  $b_1+b_2+ \dots +b_q = z_n-a_1-a_2- \dots -a_p=2n-2-2p+2=2q-2,$
  т.е. последовательность $b_1, b_2, \dots, b_q$ деревянная.

Остается доказать, что существует $s\le r$ такое что $1\le a_1 \le x_1-1$.
Для каждого $i<r$ выполнено $x_1 \ge x_i$, значит
$$z_i-i+x_1+1 \ge (z_{i+1}-(i+1)+3)-1.$$
Другими словами,

\centerline {$2=z_1-1+3$,}

\centerline {$z_1-1+x_1+1 \ge (z_2-2+3)-1$, }

\centerline {$z_2-2+x_1+1 \ge (z_3-3+3)-1$, }

\centerline {\dots, }

\centerline {$z_{r-1}-(r-1)+x_1+1 \ge (z_r-r+3)-1$,}

\centerline {$z_r-r+x_1+1=n-1$.}

Здесь последнее равенство не аналогично остальным, а следует из того, что последовательность $x_1, x_2, \dots, x_n$ деревянная и $1=x_{r+1}=\dots=x_n$.
Поскольку $2\le p\le n-1$, существует $s\le r$ такое что
$z_s -s+3\le p\le z_s-s+x_1 +1 \Leftrightarrow  1\le a_1 \le x_1-1.$ QED


{\it Доказательство части `тогда' в Теореме 2.}

Пусть   $m_1:=(n_2+n_3-n_1+1)/2, m_2:=(n_1+n_3-n_2+1)/2, m_3:=(n_1+n_2-n_3+1)/2.$
  Тогда   $m_1+m_2=n_3+1, m_1+m_3=n_2+1, m_2+m_3=n_1+1.$
По Лемме 2 существуют последовательности
  $$
p_{1 1} , p_{1 2} , \dots , p_{1 m_3}, p_{2 1} , p_{2 2} , \dots , p_{2 m_1},
 p_{3 1} , p_{3 2} , \dots , p_{3 m_2},
$$
$$
q_{1 1} , q_{1 2} , \dots , q_{1 m_2}, q_{2 1} , q_{2 3} , \dots , q_{2 m_3},
 q_{3 1} , q_{3 2} , \dots , q_{3 m_1},
$$

такие что $p_{k-1, 1}+q_{k+1, 1}=x_{k 1}$ и
упорядоченные наборы

  $(p_{k-1, 2}, p_{k-1, 3}, \dots, p_{k-1, m_{k+1}}, q_{k+1, 2}, q_{k+1, 3}, \dots, q_{k+1, m_{k-1}})$   и $(x_{k 2}, x_{k 3}, \dots, x_{k n_k})$ одинаковы с точностью до перестановок.


По Теореме 1' существуют многогранники

  $Q_1, P_1, Q_2, P_2, Q_3, P_3\subset\R^3$ такие что $Q_k \cap Q_{k+1}=\varnothing, Q_k\cap P_l=\varnothing$ если $l\neq k-1$ и

$\bullet$ $Q_k-P_{k-1}$ есть объединение $m_{k+1}$ компонент связности, $i$-тая из которых имеет $q_{k i}$ соседей


$\bullet$ $P_{k-1}-Q_k$ есть объединение $m_{k+1}$ компонент связности, $i$-тая из которых имеет $p_{k-1, i}$ соседей


$\bullet$ граница некоторой компоненты $\R^3-P_{k-1}-Q_k$ содержит компоненту $\widetilde q_k$ с $q_{k 1}$ соседями из $Q_k$ и компоненту  $\widetilde p_{k-1}$ с $p_{k-1, 1}$ соседями из $P_{k-1}$.


Для всех $1\le k \le 3$ обозначим через $f_k$ связную сумму многогранников $Q_{k+1}$ и $P_{k-1}$ и тонкой трубки, соединяющей компоненты $\widetilde q_{k+1}$ и $\widetilde p_{k-1}$ из третьего условия.

\begin{figure}[h]
  \centerline{\includegraphics[width=10cm]{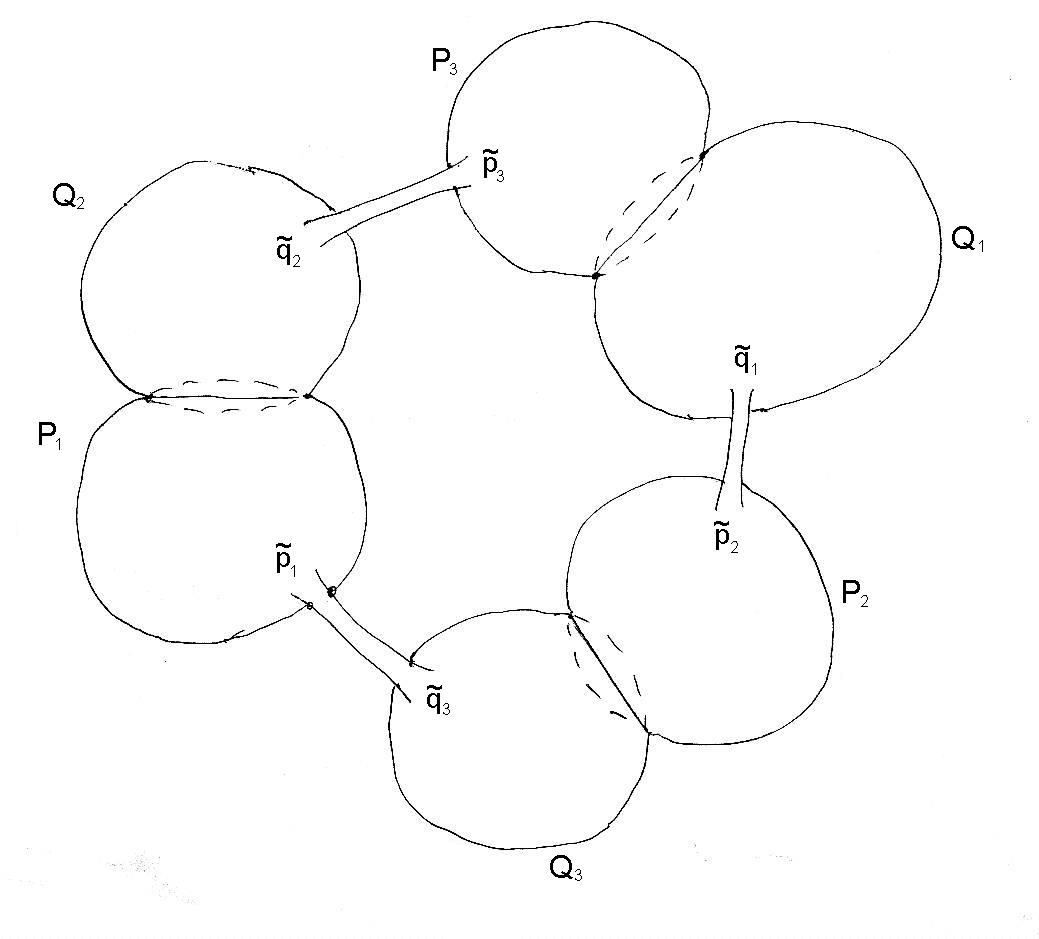}}
  \centerline{  Расположение многогранников $P_i$ и $Q_i$.}
\end{figure}

Это можно сделать без пересечений трубок.

Тогда $f_k-f_{k+1}-f_{k+2}$  будет требуемым для любого $1\le k\le 3$. QED.



\bigskip
{\bf Благодарности.}

Автор благодарен профессору С.К. Ландо за полезные замечания и обсуждения.

\bigskip
{\bf References}

[A] С. Аввакумов, Алгоритм построения двух сфер с заданным пересечением, 2012, preprint.

[D] Н. П. Долбилин, Жемчужины теории многогранников.

[H] T. Hirasa, Dissecting the torus by immersions,
Geometriae Dedicata, 145:1 (2010), 33-41.

[N] T. Nowik, Dissecting the 2-sphere by immersions,
Geometriae Dedicata 127, (2007), 37-41.

http://arxiv.org/abs/math/0612796

[R] On intersection of three embedded spheres in 3-space,
http://arxiv.org/abs/1012.0925

[W] http://ru.wikipedia.org/wiki/Многогранник